\documentclass[a4paper,12pt]{amsart}

\usepackage{fancyheadings}
\usepackage{amsmath,amssymb,amsfonts,graphicx}
\usepackage{amsmath,amssymb,amsfonts,graphicx,color,fancyhdr}
\usepackage[latin1]{inputenc}
\usepackage{psfrag}

\newtheorem{theorem}{Theorem}[section]
\newtheorem{proposition}[theorem]{Proposition}
\newtheorem{lemma}[theorem]{Lemma}
\newtheorem{definition}[theorem]{Definition}

\newtheorem{remark}[theorem]{Remark}

\usepackage[colorlinks=true, pdfstartview=FitV, linkcolor=blue, citecolor=blue, urlcolor=blue]{hyperref}
\usepackage{color}
\usepackage{psfrag,caption}

\def\be#1 {\begin{equation} \label{#1}}
\newcommand{\ee}{\end{equation}}

\def\sqw{\hbox{\rlap{\leavevmode\raise.3ex\hbox{$\sqcap$}}$%
\sqcup$}}
\def\findem{\ifmmode\sqw\else{\ifhmode\unskip\fi\nobreak\hfil
\penalty50\hskip1em\null\nobreak\hfil\sqw
\parfillskip=0pt\finalhyphendemerits=0\endgraf}\fi}

\newcommand{\R}{{\mathbb {R}}}

\newcommand{\N}{{\mathbb N}}
\newcommand{\Z}{{\mathbb Z}}
\newcommand{\C}{{\mathbb C}}

\newcommand{\T}{{\mathbb T}}

\begin{document}

\title
[A note on bi-linear multipliers]
{\bf A note on bi-linear multipliers}

\date{\today}


\author{Saurabh Shrivastava}
\address{Department of mathematics, Indian Institute of Science Education and Research Bhopal, Indore By-pass road Bhauri, Bhopal-462066, India}
\email{saurabhk@iiserb.ac.in}
\begin{abstract} In this paper we prove that if $\chi_{_E}(\xi-\eta)-$ the indicator function of measurable set $E\subseteq \R^d,$ is a bi-linear multiplier symbol for exponents $p,q,r$ satisfying the H\"{o}lder's condition $\frac{1}{p}+\frac{1}{q}=\frac{1}{r}$ and exactly one of $p,q,$ or $r'=\frac{r}{r-1}$ is less than $2,$ then $E$ is equivalent to an open subset of $\R^d.$   
\end{abstract}
\maketitle


\section{Introduction and statement of results}
The remarkable work of M.~Lacey and C.~Thiele
\cite{lt1}, \cite{lt2} on boundedness of the bi-linear Hilbert transform motivated a lot of research in the area of Euclidean harmonic analysis. For $f,g$ in $\mathcal{S}(\R)-$the Schwartz class on $\R$, the bi-linear Hilbert transform is defined
by \begin{eqnarray*}\label{bht1}
H(f, g)(x) = p.v. \int_{\R} f(x-y)g(x+y)\frac{dy}{y}.
\end{eqnarray*}
Or equivalently, 
\begin{eqnarray*}\label{bht2}
H(f,g)(x)= -i \int_{\R}\int_{\R}\hat{f}(\xi)\hat{g}(\eta) \it{sgn}(\xi-\eta)e^{2\pi i x(\xi+\eta)}d\xi d\eta,
\end{eqnarray*}
where $\hat{}$ denotes the Fourier transform and $$\it{sgn}(\xi)=\left\{
  \begin{array}{ll}
    1, & \hbox{$\xi>0$} \\
    0, & \hbox{$\xi=0$}\\
   -1, & \hbox{$\xi<0$.}
  \end{array}
\right.$$

We would like to remark that the bi-linear Hilbert transform is invariant under the operations of simultaneous translation, dilation, and modulation. The modulation in-variance is a subtle property shared by the bi-linear Hilbert transform and poses additional difficulties while proving suitable $L^p-$estimates for the operator. It is not difficult to convince ourselves that the classical approach of Littlewood-Paley decomposition, which is a useful technique to handle singular integral operators, is not quite helpful to deal with operators having modulation in-variance property. In papers~\cite{lt1} and \cite{lt2} M.~Lacey and C.~Thiele developed techniques very useful to handle such operators. These techniques are commonly referred to as the time-frequency techniques. In seminal papers they proved the following $L^p-$estimates for the bi-linear Hilbert transform. 
\begin{theorem}\cite{lt1,lt2}\label{lt}
Let $1<p,q\leq \infty$ and $\frac{2}{3}<r <\infty$ be such that $\frac{1}{p}+\frac{1}{q}=\frac{1}{r}$. Then for all functions $f,g \in \mathcal{S}(\R)$, there exists a constant $C>0$ such that
$$\|H(f,g)\|_{L^{r}(\R)}\leq C \|f\|_{L^{p}(\R)}\|g\|_{L^{q}(\R)}.$$
\end{theorem}
In this paper we are interested in studying bi-linear multiplier operators having modulation in-variance property. The bi-linear multiplier operators in general are defined as follows: Let $m(\xi-\eta)$ be a bounded measurable function on $\R^d$ and $(p,q,r),~0<p,q,r\leq \infty$ be a triplet of exponents.  Consider the bi-linear operator $T_m$ initially defined for functions $f$ and $g$ in a suitable dense class by  
\begin{eqnarray}\label{def1}
T_{m}(f,g)(x)= \int_{\R^d}\int_{\R^d}\hat{f}(\xi)\hat{g}(\eta) m(\xi-\eta)e^{2\pi i x\cdot(\xi+\eta)}d\xi d\eta.
\end{eqnarray}

We say that $T_m$ is a bi-linear multiplier operator for the triplet $(p,q,r)$ if $T_m$ extends to a bounded operator 
from $L^p(\R^d)\times L^q(\R^d)$ into $L^r(\R^d),$ i.e. there exists a constant $C>0,$ independent of functions $f$ and $g,$ such that  
\begin{eqnarray*}\label{bm}
\|T_m(f,g)\|_{L^r(\R^d)}\leq C \|f\|_{L^p(\R^d)} \|g\|_{L^q(\R^d)}.
\end{eqnarray*}

The bounded function $m$ is said to be a bi-linear multiplier symbol for the triplet $(p,q,r)$ if the corresponding operator $T_m$ is a bi-linear multiplier operator for $(p,q,r).$ 

We denote by $\mathcal{M}_{p,q}^{r}(\R^d)$ the space of all bi-linear multiplier symbols for the triplet $(p,q,r).$ Further,  the norm of $m\in \mathcal{M}_{p,q}^{r}(\R^d)$ is defined to be the norm of the corresponding bi-linear multiplier operator $T_m$ from $L^p(\R^d)\times L^q(\R^d)$ into $L^r(\R^d),$ i.e.
\begin{eqnarray*}
\|m\|_{\mathcal{M}_{p,q}^{r}(\R^d)}=\|T_m\|_{L^p(\R^d)\times L^q(\R^d)\rightarrow L^r(\R^d)}.
\end{eqnarray*}

The bi-linear multiplier symbols on the Torus group $\T^d$ and discrete group $\Z^d$ are defined similarly. The space of bi-linear multiplier symbols on $\T^d$ and $\Z^d$ will be denoted by $\mathcal{M}_{p,q}^{r}(\T^d)$ and $\mathcal{M}_{p,q}^{r}(\Z^d)$ respectively. 
\begin{remark}Unless specified otherwise, we shall always assume that exponents $p,q,r$ satisfy $0< p,q,r\leq \infty$ and the H\"{o}lder's condition $\frac{1}{p}+\frac{1}{q}=\frac{1}{r}.$
\end{remark}
Now we describe some important properties of bi-linear multipliers. Some of these properties will be used later in the paper. 
\begin{proposition}\cite{g}\label{gra} Let $p,q,r$ be exponents satisfying the H\"{o}lder's condition. Then bi-linear multiplier symbols satisfy the following properties:
\begin{enumerate}
\item If $c\in \C$ and $m, m_1,$ and $m_2$ are in $\mathcal{M}_{p,q}^{r}(\R^d),$ then so are $c m$ and $m_1+m_2.$ Moreover, 
$\|cm\|_{\mathcal{M}_{p,q}^{r}(\R^d)}=|c|~\|m\|_{\mathcal{M}_{p,q}^{r}(\R^d)} $ and $\|m_1+m_2\|_{\mathcal{M}_{p,q}^{r}(\R^d)}\leq C(\|m_1\|_{\mathcal{M}_{p,q}^{r}(\R^d)} +\|m_2\|_{\mathcal{M}_{p,q}^{r}(\R^d)}).$
\item If $m\in \mathcal{M}_{p,q}^{r}(\R^d)$ and $\eta\in \R^d,$ then $\tau_{\eta}m(.)=m(.-\eta)$ is in $\mathcal{M}_{p,q}^{r}(\R^d)$ with same norm as $m.$
\item If $m\in \mathcal{M}_{p,q}^{r}(\R^d)$ and $\lambda>0,$ then $m_{\lambda}(.)=m(\lambda.)$ is in $\mathcal{M}_{p,q}^{r}(\R^d)$ with same norm as $m.$
\item If $m\in \mathcal{M}_{p,q}^{r}(\R^d)$ and $A=(a_{i,j})_{d\times d}$ is an orthogonal matrix acting on $\R^d,$ then $m(A.)$ is in $\mathcal{M}_{p,q}^{r}(\R^d)$ with same norm as $m.$
\item If $m\in \mathcal{M}_{p,q}^{r}(\R^d)$ and $h\in L^1(\R^d),$ then the convolution $m\ast h\in \mathcal{M}_{p,q}^{r}(\R^d),$ provided $r\geq 1.$ Moreover, we have $\|m\ast h\|_{\mathcal{M}_{p,q}^{r}(\R^d)}\leq \|h\|_{L^1(R^d)}~\|m\|_{\mathcal{M}_{p,q}^{r}(\R^d)}.$
\item If $m_j$ is a sequence of bi-linear multiplier symbols in $\mathcal{M}_{p,q}^{r}(\R^d)$ such that the norm 
$\|m_j\|_{\mathcal{M}_{p,q}^{r}(\R^d)}\leq C$ for some fixed constant $C>0$ uniformly for all $j=1, 2, \dots.$ Further assume that $m_j$ are uniformly bounded by a bounded measurable function on $\R^d$, then if $m_j\rightarrow m$ a.e. as $j\rightarrow \infty,$ we have 
that $m\in \mathcal{M}_{p,q}^{r}(\R^d)$ and $\|m\|_{\mathcal{M}_{p,q}^{r}(\R^d)}\leq C.$
\end{enumerate}
\end{proposition}
\begin{remark}
We would like to remark that corresponding properties hold true for bi-linear multiplier symbols on $\T^d$ and $\Z^d.$ 
\end{remark}

\subsection{Statement of the result}
In this paper we study some particular type of bi-linear multiplier symbols, namely those which are indicator functions of measurable sets. In general, there is no effective method to decide that the indicator function of a measurable set is a bi-linear multiplier symbol for some exponents. An important example in this direction is given by $\chi_{_I}(\xi-\eta),$ for an interval $~I\subset \R.$ This is a consequence of boundedness of the bi-linear Hilbert transform that $\chi_{_I}(\xi-\eta)$ is a bi-linear multiplier symbol for all exponents satisfying the conditions of Theorem~\ref{lt}. 

In the current paper, we study structural properties (in the sense of measure theory) of sets whose indicator functions give rise to bi-linear multiplier symbols. The motivation for this paper comes from the beautiful work of V.~Lebedev and A.~Olevski\u{i}~\cite{lo} on the classical Fouirer multipliers. We will prove an analogue of their result in the context of bi-linear multiplier operators. In order to describe  V.~Lebedev and A.~Olevski\u{i}'s result, we need the following definition. 
\begin{definition}
We say that measurable sets $E$ and $E'$ are equivalent if the symmetric difference $E\Delta E'$ has Lebesgue measure zero. 
\end{definition}
\begin{theorem}\label{leb}\cite{lo} Let $E\subseteq \R^d$ be a measurable set and $p\neq 2.$ If $\chi_{_E}-$the indicator function of $E,$ is an $L^p-$multiplier, i.e., the linear operator 
$f\rightarrow (\chi_{_E}\hat f\check{)},~f\in L^2(\R^d)\cap L^p(\R^d)$ extends boundedly from $L^p(\R^d)$ into itself, then $E$ is equivalent to an open set in $\R^d.$  
\end{theorem} 
This theorem tells us that the structure of set whose indicator function is an $L^p-$multiplier, $p\neq 2,$ cannot be very complicated in the sense of measure theory. As an immediate consequence of this, we see that the indicator function of a nowhere dense set of positive Lebesgue measure is never an $L^p-$multiplier for $p\neq 2.$

As mentioned previously, in this paper our aim is to prove an analogue of Theorem~\ref{leb} in the context of bi-linear multipliers. In particular, we prove the following result:  
\begin{theorem}\label{result}
Let $E$ be a non-empty measurable subset of $\R^d$ and $p,q,r$ be exponents such that $1\leq p,q,r\leq \infty,~\frac{1}{p}+\frac{1}{q}=\frac{1}{r}$ and exactly one of $p,q,$ or $r'=\frac{r}{r-1}$ is less than $2.$ Suppose that $\chi_{_E}(\xi-\eta)$ is a bi-linear multiplier symbol for the triplet $(p,q,r)$, then $E$ is equivalent to an open subset of $\R^d.$ 
\end{theorem}

\begin{remark} 
In the case of classical Fourier multipliers $p=2$ plays a special role by virtue of the Plancherel theorem.  In a sharp contrast to this, in the theory of bi-linear mutipliers there is absolutely no easy way by which one can test that a given bounded function $m(\xi-\eta)$ is a bi-linear multiplier symbol for a given triplet of exponents. The range of exponents covered by Theorem~\ref{result} falls in the complement of what is commonly known as local $L^2-$range of exponents. The local $L^2-$range consists of exponents satisfying $2\leq p,q,r'\leq 2.$ Generally, it is believed that in this range of exponents most of the bi-linear multiplier operators are well behaved as far as the boundedness is concerned. 
\end{remark}
\section{Basic results and proof of Theorem~\ref{result}}
In this section first we provide some basic definitions and results which would be required to complete the proof of Theorem ~\ref{result}.
\begin{definition}Let $E$ be a measurable set of $\R^d.$ Then we say that a point $x\in \R^d$ is a density point for $E$ if 
\begin{eqnarray}
\lim\limits_{t\rightarrow 0}\frac{|B(x,t)\cap E|}{|B(x,t)|}=1,
\end{eqnarray}
where $B(x,t)$ denotes the Euclidean ball of radius $t>0$ centered at $x\in \R^n$ and $|.|$ denotes the Lebesgue measure of a set.
 
The set of all density points for the set $E$ is denoted by $E^d.$ The set $\partial_e E=\overline{E^d}\cap \overline{(E^c)^d}$ is referred to as the  essential boundary of $E,$ where $E^c$ denotes the complement of $E.$ 
\end{definition}
\begin{lemma}\label{l0} If $E\subseteq \R^d$ is a measurable set, then $E$ and $E^c$ are both equivalent to open sets if and only if $\partial_e E-$the essential boundary of $E,$ has Lebesgue measure zero. 
\end{lemma}
This lemma is easy to verify and hence its proof is not included here. Next, we describe an important lemma proved by V.~Lebedev and A.~Olevski\u{i} in~\cite{lo}. 
\begin{lemma}\label{leb1}\cite{lo} Let $E\subseteq \R^d$ be a measurable set such that $|\partial_e E|>0.$ Then for every $N\in \N$ and for every subset $A\subseteq A_N=\{1,2,\dots,N\},$ there exist $x_0, h\in \R^d$ such that the arithmetic progression 
$$x_n=x_0+nh,~n\in A_N$$
satisfies the following conditions 
\begin{eqnarray}\label{l1}
x_n\in E^d,~~\text{if}~~n\in A~~\text{and}~~x_n\in (E^c)^d,~\text{if}~~n\in A_N\setminus A.
\end{eqnarray}
\end{lemma} 
This lemma plays a crucial role in the proof of Theorem~\ref{leb}. Since, we have exploited the methodology of~\cite{lo} in order to prove the main result of this paper, Lemma~\ref{leb1} is a key tool for the current paper as well. 

Finally, we would require de Leeuw's type transference result for bi-linear multipliers. This would allow us to restrict bi-linear multiplier symbols defined on $\R^d$ to the discrete group $\Z^d.$ As a consequence of this, we will have to deal with some bi-linear multiplier operators on $\T^d.$ In the context of this paper, it turns out that working with operators on $\T^d$ is simpler compared to working with the original operator defined on $\R^d.$ This approach helps us in proving some good estimates on bi-linear operators under consideration and eventually helps us completing the proof of our main result.

The following is a bi-linear analogue of the celebrated transference result proved by de Leeuw~\cite{de} for the classical Fourier multipliers.  
\begin{theorem}\label{ob}\cite{oscar}
Let $m(\xi-\eta) \in \mathcal{M}_{p,q}^{r}(\R^d)$ be a regulated function. Then $m\restriction_{\Z^d}$ belongs to $\mathcal{M}_{p,q}^{r}(\T^d)$ with norm bounded by a constant multiple of $\|m \|_{ \mathcal{M}_{p,q}^{r}(\R^d)}.$
\end{theorem}

We are now in a position to prove the main result of this paper. 
\subsection{Proof of Theorem~\ref{result}}
Notice that the same property holds for the complement $E^c,$ i.e., the indicator function $\chi_{_{E^c}}$ is a bi-linear multiplier symbol for the triplet $(p,q,r).$ Therefore, in the view of Lemma~\ref{l0}, it is enough to 
show that the essential boundary $\partial_e E$ has Lebesgue measure zero. 

The proof is given by contradiction. Therefore, suppose on the contrary that $|\partial_e E|>0.$ 
Let $N\in \N$ be fixed and $\{\epsilon_n\}_{n=1}^N$ be a random sequence of $0$ and $1.$ Then as an application of Lemma~\ref{leb1} we know that there exists $x_0, h\in \R^d$ such that 
$$x_n=x_0+nh \in \left\{
  \begin{array}{ll}
    E^d, & \hbox{$\epsilon_n=1$} \\
    (E^c)^d, & \hbox{$\epsilon_n=0$.}
  \end{array}
\right.$$
For $t>0,$ we consider the following function 
\begin{eqnarray*}m_t(x)&=&\frac{1}{|B(x,t)|}\int_{B(x,t)}\chi_{_E}(y)dy\\
&=& \frac{|B(x,t)\cap E|}{|B(x,t)|}.
\end{eqnarray*}
Observe that $m_t=\chi_{_E}\ast A_t,$ where $A_t(.)=\frac{1}{|B(0,t)|}\chi_{_{B(0,t)}}(.).$ Clearly, for all $t$ the norm $\|A_t\|_{L^{1}(\R^d)}=1.$  Since $r\geq 1,$ by Proposision~\ref{gra} we obtain that 
$m_t\in \mathcal M_{p,q}^r(\R^d).$ Moreover 
$$\|m_t\|_{\mathcal M_{p,q}^r(\R^d)}\leq \|\chi_{_E}\|_{\mathcal M_{p,q}^r(\R^d)},~\forall t>0.$$ 

Let $e_1=(1,0,\dots,0)$ be the unit vector of the $x_1-$axes and $\psi$ be an affine mapping of $\R^d$ which maps the vector $ne_1$ to $x_n.$ At this point we invoke Proposition~\ref{gra} once again and obtain that the composition $m_t\circ \psi\in {\mathcal M_{p,q}^r}(\R^d)$ with uniform norm with respect to the parameter $t.$ 

Recall the definition of density points and notice that $m_t(x_n)\rightarrow \epsilon_n$ as $t\rightarrow 0,$ which 
is the same as $m_t\circ\psi(ne_1)\rightarrow \epsilon_n$ as $t\rightarrow 0.$ 

Now we identify the set $\{ne_1:n\in \Z\}$ with the discrete group $\Z$ and invoke de Leeuw's type transference theorem for bi-linear multipliers, namely Theorem~\ref{ob} from~\cite{oscar}. This yields that $m_t\circ \psi\restriction_{\Z}-$the restriction of $m_t\circ \psi$ to $\Z,$ belongs to $\mathcal M_{p,q}^r(\T)$ and  
$$\|m_t\circ \psi\restriction_{\Z}\|_{\mathcal M_{p,q}^r(\T)}\leq C\|m_t\|_{\mathcal M_{p,q}^r(\R^d)}~\forall t>0.$$ 
Since the above estimate is uniform with respect to $t>0,$ a standard limiting argument (see Proposition~\ref{gra}) allows us to conclude that the sequence $\{\epsilon_n\}_{n=1}^N$ is a bi-linear multiplier symbol for the triplet $(p,q,r),$ i.e., $\{\epsilon_n\}\in {\mathcal M_{p,q}^r}(\T).$ Moreover the norm is independent of $N\in \N.$ 

We shall show that this leads to a contradiction.  

First, we will prove that exponent $p$ or $q$ cannot be smaller than $2.$ Since $\{\epsilon_n\}\in {\mathcal M_{p,q}^r}(\T)$ with multiplier norm uniformly bounded in $N\in \N,$ we use Khintchine's inequality to deduce that the bi-linear Littlewood-Paley operator 
$$S:(P,Q)\rightarrow \left(\sum_n |S_n(P,Q)|^2\right)^{\frac{1}{2}}$$ 
is bounded from $L^p(\T)\times L^q(\T)$ into $L^r(\T),$ where $P$ and $Q$ are trigonometric polynomials and $S_n$ is the bi-linear multiplier operator on $\T$ with corresponding bi-linear multiplier symbol $\chi_{_{\{n\}}}.$ But, we already know that $p,q\geq 2$ is a necessary condition for the boundedness of the bi-linear Littlewood-Paley operator $S$ (see the paper by P.~Mohanty and S.~Shrivastava~\cite{ms2} for a proof of this assertion) and hence we get a contradiction if either of $p$ and $q$ is less than $2.$ 

Next, we assume that $r'<2,$ and show that this assumption also gives a contradiction. Let $p,q,r$ be exponents such that $p,q,r>2.$  

Since the sequence $\{\epsilon_n\}\in {\mathcal M_{p,q}^r}(\T)$ with multiplier norm uniformly bounded in $N\in \N,$ we have that  
\begin{equation}\label{lin}
\|S_{\Gamma}(P,Q)\|_{L^{r}(\T)} \leq C \|P\|_{L^{p}(\T)} \|Q\|_{L^{q}(\T)},
\end{equation}
where $S_{\Gamma}$ is the bi-linear multiplier operator on $\T$ with corresponding bi-linear multiplier symbol $\{\epsilon_n\}.$ 

Take $Q\equiv 1$ on $\T$ and consider
\begin{eqnarray*}
S_\Gamma(P,Q)(x)&=&\sum_{m,l\in \Z} \hat P(m)\hat Q(l) \epsilon_{l-m} e^{2 \pi i (m+l)x}\\
&=&\sum_{m\in \Z} \hat P(m)\epsilon_{-m} e^{2 \pi i mx}.
\end{eqnarray*}
Therefore, we obtain that the linear operator $T:f\rightarrow  \sum_{n\in \Z} \epsilon_{n} \hat P(n) e^{2 \pi i n\cdot}$ is bounded from $L^p(\T)$ into $L^r(\T),$ i.e., we have 
\begin{eqnarray}\label{lin1}
\|T P\|_{L^r(\T)}\leq C \|P\|_{L^p(\T)}.
\end{eqnarray} 
This leads to a contradiction as we know that the inequality~(\ref{lin1}) does not hold true for all choices of $\{\epsilon_n\}_{n=1}^N,~N\in \N$ as $r>2.$\\
This completes the proof of Theorem~\ref{result}. 

\end{document}